\documentclass{article}

\usepackage{amssymb, amsmath, amsthm, amsfonts, amscd, epsf, graphicx}
\usepackage[dvips,usenames]{color}
\usepackage[american]{babel}
\usepackage[latin1]{inputenc}
\usepackage[T1]{fontenc}
\usepackage{ae, aecompl}
\usepackage[mathscr]{eucal}
\usepackage{graphicx}  
\usepackage{epsfig}
\usepackage{float}

\newtheorem{theorem}{Theorem}
\newtheorem{lemma}{Lemma}
\newtheorem{corollary}{Corollary}

\newtheorem{remark}{Remark}

\title{Computational  aspects of   hyperelliptic curves}

\author{T. Shaska\\
Department of Mathematics, \\
University of California at  Irvine \\
E-mail: tshaska@math.uci.edu }
\date{}

\begin{document}

\maketitle



\def\Z{\mathbb Z}
\def\bZ{\mathbb Z}
\def\Q{\mathbb Q}
\def\C{\mathbb C}
\def\bP{\mathbb P}
\def\cB{\mathcal B}
\def\cA{\mathcal A}
\def\L{\mathcal L}
\def\cR{\mathcal R}
\def\H{\mathcal H}
\def\M{\mathcal M}
\def\N{\mathcal N}
\def\w{\widetilde}
\def\l{\lambda}
\def\s{\sigma}
\def\G{\Gamma}
\def\a{\alpha}
\def\b{\beta}
\def\p{\mathfrak p}
\def\P{\mathcal P}
\def\e{\varepsilon}
\def\iso{\equiv}
\def\sem{{\rtimes}}

\def\bG{\overline G}
\def\g{\gamma}
\def\bg{\bar \gamma}
\def\u{\mathfrak u}
\def\k{\bar k}
\def\iso{{\, \cong\, }}
\def\nor{{\,  \vartriangleleft \, }}
\def\<{\langle}
\def\>{\rangle}
\def\emb{\hookrightarrow }


%
\begin{abstract}
We introduce a new approach of computing the automorphism group and the field of moduli of  points $\p=[C]$ in the
moduli space of hyperelliptic curves $\H_g$. Further, we show that for every moduli point $\p \in \H_g(L)$ such
that the reduced automorphism group of $\p$ has at least two involutions, there exists a representative $C$ of the
isomorphism  class $\p$ which is defined over $L$.
\end{abstract}

\section{Introduction}
The  purpose of  this  note is  to  introduce some  new techniques  of computing the automorphism group and  the
field of moduli of genus $g$ hyperelliptic curves.  Former results by many authors  have focused on hyperelliptic
curves  of  small  genus,  see
 \cite{Me}, \cite{CQ}, \cite{MS}, \cite{Sh2}, \cite{SV},
et. al.  We aim to find a method which would work for any genus.

Let  $C$  denote a  genus  $g$  hyperelliptic  curve defined  over  an algebraically closed field $k$  of
characteristic zero and $G:=Aut(C)$ its automorphism  group.      We denote by $\H_g$  the moduli space of genus
$g$  hyperelliptic curves and by  $\L_g$ the locus  in $\H_g$ of hyperelliptic   curves   with   extra
involutions.   $\L_g$   is   a $g$-dimensional rational variety,  see \cite{GS2}. Equation  \ref{eq} gives a
normal form  for curves  in $\L_g$.   This normal  form  depends on parameters $a_1, \dots , a_g \in k$, such that
the discriminant of the right  side $\Delta(a_1, \dots  , a_g)  \neq 0$.   Dihedral invariants $(u_1, \dots ,
u_g)$ were introduced by Gutierrez  and this author in \cite{GS2}. The  tuples  $\u=(\u_1,  \dots  , u_g)$ (such
that $\Delta_{\u}\neq 0$) are in one-to-one correspondence with isomorphism classes of genus $g$ hyperelliptic
curves with automorphism group the Klein 4-group.  Thus,  dihedral invariants $u_1, \dots ,  u_g$ yield a
birational  parameterization  of  the locus  $\L_g$.   Computationally these invariants give an efficient way of
determining a generic point of  the  moduli  space  $\L_g$.   Normally, this  is  accomplished  by invariants of
$GL_2(k)$ acting on  the space of binary forms of degree $2g+2$.   These $GL_2(k)$-invariants  are not  known for
$g  \geq 3$. However, dihedral invariants are explicitly defined for all genera.

The  full automorphism  groups   of   hyperelliptic  curves are determined in \cite{Bu} and  \cite{BS}. Most of
these groups have non-hyperelliptic involutions (i.e., the corresponding curve is in $\L_g$). For each group  $G$
that occurs as full  automorphism  group  of  genus  $g$ curves  one  determines  the $G$-locus  in $\L_g$  in
terms  of the  dihedral invariants.   Given a genus $g$ curve $C$  we first determine if $C\in \L_g$. Then we
compute its dihedral invariants and determine the locus $\L_G$ that they satisfy. This determines  $Aut(C)$.
Present algorithms of  computing the automorphism group of a hyperelliptic curve $Y^2=F(X)$ are based on computing
the roots of $F(X)$ and then finding fractional  linear  transformations  that  permute  these  roots.   The
algorithm we propose requires   only determining the normal form of $C$ (i.e., Eq.~\ref{eq}). This requires
solving a system of $g$-equations and four unknowns. For curves which have at least two involutions in their
reduced automorphism group we find a nice condition on the dihedral invariants.

For $C\not \in \L_g$ similar methods can be used. If  $|Aut(C)|  >2$ and $C \not \in \L_g$,  then $C$ has an
automorphism of order $N$, where $N$ is as in Lemma \ref{lemma3.3}. For small genus these curves can be classified
by ad-hoc methods. In general one needs to find invariants of such spaces for all $N>2$ and implement similar
methods as above. We intend this as the  object of further research.

In section 4,  we introduce  how to compute the field of moduli of genus $g$ hyperelliptic curves with
automorphism group of order $> 4$. Let $\M_g$  (resp., $\H_g$)  be the moduli  space of  algebraic curves (resp.,
hyperelliptic curves)  of genus  $g$  defined over  $k$  and $L$  a subfield of  $k$.  It is well  known that
$\M_g$ (resp.,  $\H_g$) is a $3g-3 $  (resp., $2g-1$) dimensional variety.   If $C$ is  a genus $g$ curve defined
over $L$, then clearly $[C]\in  \M_g(L)$.  However, the converse  is not  true.  In  other words,  the moduli
space  $\M_g$ of algebraic curves of genus $g$ is a coarse moduli space.  The answer is not   obvious  if   we
restrict ourselves  to   the  singular   points  of   $ \M_g$.  Singular points  of  $\M_g$ (resp.,  $\H_g$)
correspond  to isomorphism  classes  of  curves  with nontrivial  automorphism  groups (resp.,  automorphism
groups  of  order  $>  2$  ).   In  general,  we conjecture that for a singular  point $\p \in \M_g(L)$ (resp.,
$\p \in \H_g(L)$)  there  is  always  a  curve  $C$  defined  over  $L$  which correspond to $\p$.  We focus on
$\H_g$.  A point $\p=[C] \in \H_g$ is given by the $g$-tuple of  dihedral invariants. We denote by $Aut(\p)$ the
automorphism group  of any  representative  $C $  of $\p$.   More precisely, for hyperelliptic curves we
conjecture the following:

\smallskip

\noindent {\bf  Conjecture 1:} {\it Let $\p  \in \H_g  (L)$ such  that $|Aut(\p)|  > 2$. There exists a
representative $C$ of the isomorphism  class $\p$ which is defined over $L$. }

\smallskip

\noindent In this  paper we show  how dihedral invariants  can be used  to prove some special cases  of this
conjecture. A detailed  discussion on this problem is intended in \cite{Sh3}.  The condition $|Aut(\p)| > 2$ of
the above conjecture  can not be  dropped.  Determining exactly  the points $\p \in  \H_g\setminus \L_g$  where
such rational  model $C$  does not exist is  still an  open problem.  For  $g=2$ Mestre (1991)  found an algorithm
which determines  such points.   It is  based  on classical invariants of binary sextics.

\medskip

\noindent  {\bf  Notation:}  Throughout  this  paper  $k$  denotes  an algebraically  closed field  of
characteristic zero,  $g$ an  integer $\geq 2$, and $C$ a  hyperelliptic curve of genus $g$.  $\M_g$ (resp.,
$\H_g$) is  the moduli space  of curves (resp.,  hyperelliptic curves) defined  over  $k$.  Further,  $V_4$
denotes  the  Klein 4-group  and $D_{2n}$  (resp., $\Z_n$)  the  dihedral group  of  order $2n$  (resp., cyclic
group of order $n$).

\section{Dihedral invariants of hyperelliptic curves}

Let $k$  be an algebraically  closed field of characteristic  zero and $C$  be  a  genus  $g$  hyperelliptic
curve  given  by  the  equation $Y^2=F(X)$, where $\deg(F)=2g+2$. Denote  the function field of $C$ by
$K:=k(X,Y)$.  Then, $k(X)$ is the  unique degree 2 genus zero subfield of  $K$.   We  identify  the  places  of
$k(X)$  with  the  points  of $\bP^1= k \cup \{\infty\}$  in the  natural  way (the  place $X=\a$  gets identified
with  the point $\a \in \bP^1$).  Then, $K$ is  a quadratic extension field of $k(X)$  ramified exactly at
$n=2g+2$ places $\a_1, \dots , \a_n$  of $k(X)$.  The corresponding places  of $K$ are called the {\it Weierstrass
points} of $K$.  Let $\P:=\{\a_1, \dots , \a_n \}$. Thus, $K=k(X,Y)$, where
%
%

\begin{equation}\label{eq1}
Y^2=\prod_{{\a\in \P}, \, \, { \a \neq \infty} }  (X-\a).
\end{equation}

Let $G=Aut(K/k)$.  Since $k(X)$  is the only  genus 0 subfield  of   degree  2  of   $K$,  then  $G$  fixes
$k(X)$.  Thus, $G_0:=Gal(K/k(X))=\< z_0 \>$, with $z_0^2=1$, is central in $G$. We call {\bf the  reduced
automorphism group}  of $K$ the  group $\bG:=G/G_0$. Then, $\bG$  is naturally isomorphic to the  subgroup of
$Aut(k(X)/k)$ induced  by  $G$.  We  have a  natural  isomorphism  $\Gamma:=PGL_2(k) \to Aut(k(X)/k)$. The action
of  $\Gamma$ on the places of  $k(X)$ corresponds under the above  identification to  the usual  action on
$\bP^1$  by fractional linear  transformations $t \mapsto  \frac {at+b}  {ct+d}$. Further, $G$ permutes  $\a_1,
\dots  ,  \a_n$. This  yields  an  embedding  $\bG \emb S_n$.

Because  $K$ is  the  unique  degree 2  extension  of $k(X)$  ramified exactly  at  $\a_1$,  \dots  ,  $\a_n$,
each  automorphism  of  $k(X)$ permuting these  $n$ places extends  to an automorphism of  $K$. Thus, $\bG$ is the
stabilizer in  $Aut(k(X)/k)$ of the set $\P$. Hence under the isomorphism $\Gamma \mapsto Aut(k(X)/k)$, $\bG$
corresponds to the stabilizer $\Gamma_\P$ in $\Gamma$ of the $n$-set $\P$.

An extra involution of $K$ is  an involution in $G$ which is different from  $z_0$ (the  hyperelliptic
involution).  If  $z_1$  is an  extra involution and  $z_0$ the hyperelliptic one, then  $z_2:=z_0\, z_1$ is
another extra involution.  So  the extra involutions come naturally in pairs.  Suppose $z_1$ is an  extra
involution of $K$. Let $z_2 :=z_1\, z_0$,  where $z_0$  is the  hyperelliptic involution.  Then $K=k(X,Y)$ with
equation
\begin{equation}\label{eq}
Y^2=X^{2g+2} + a_{g} X^{2g}+ \dots + a_1 X^2 +1
\end{equation}
see  \cite{GS2}. The  dihedral  group $H:=D_{2g+2}=\<\tau_1, \tau_2 \>$  acts on  $k(a_1, \dots , a_g)$ as
follows:
\begin{eqnarray*}
\tau_1: & \quad a_i & \to \e^{2i} a_i, \qquad  for \quad i=1, \dots , g \\
\tau_2:  & \quad a_i  &\to a_{g+1-i}, \qquad for \quad  i=1, \dots , [\frac {g+1} 2]
\end{eqnarray*}
The fixed  field $k(a_1,  \dots ,  a_g)^H$ is  the same  as the
 function field of the variety $\L_g$. The invariants of
 such action are
\begin{equation}
u_i:=  a_1^{g-i+1} \, a_i \, + \, a_g^{g-i+1} \, a_{g-i+1}, \quad for \quad 1 \leq i \leq g
\end{equation}
and are called {\bf  dihedral invariants} for the genus  $g$ and the tuple
$$\u:=(u_1,  \dots ,  u_g)$$  is  called the  {\bf  tuple of  dihedral
invariants}, see \cite{GS2} for details.

It is easily seen that $\u=0$ if and only if $a_1=a_g=0$. In this case replacing $a_1, a_g$ by $a_2, a_{g-1}$ in
the formula above would give new invariants.  In  \cite{GS2} it is shown that  $k(\L_g)=k(u_1, \dots , u_g)$.
The  $(2g+2)$-degree  field extension $k(a_1, \dots , a_g) /  k(u_1, \dots , u_g)$ has equation
\begin{equation}
2^{g+1}\, a_g^{2g+2} - 2^{g+1}\, u_1 \,  a_g^{g+1} + u_g^{g+1}=0
\end{equation}
and the map
\begin{eqnarray*}
  \Phi:   &    k\setminus \{\Delta \neq 0\}   & \to \L_g\\
          &   (a_1, \dots , a_g)         &   \to (u_1, \dots , u_g)
\end{eqnarray*}
has  Jacobian  zero  exactly on points which correspond to curves $C\in \L_g$ such that $V_4 \emb \bG$.

\section{Automorphism groups}
In this section we  suggest an algorithm for computing the full automorphism group of hyperelliptic curves. Let
$C$ be a genus $g$ hyperelliptic curve with equation $Y^2=F(X)$ where $\deg (F)=2g+2$. Existing algorithms are
based on finding all automorphisms of $C$. Instead, we search  for only one  automorphism (non-hyperelliptic) of
$C$ of  order $N$. Most of the  time $N=2$ is enough since the majority of groups of order $> 2$ that occur as
full automorphism groups have non-hyperelliptic involutions. It is well known that the order of a non-trivial
automorphism of a hyperelliptic curve is $2 \leq N \leq 2(2g+1)$, where $2(2g+1)$ is known as the Wiman's bound.

If an automorphism of order $N=2$ exists then $C \in \L_g$ and we use dihedral invariants to determine the
automorphism group. We illustrate with curves of small genus.

The case  $g=2$ has been studied in \cite{SV}. Every point in $\M_2$ is a triple $(i_1, i_2, i_3)$ of absolute
invariants.  We state the results of \cite{SV}  without proofs.
\begin{lemma} \label{u_v}
Let $C$ be a genus 2 curve such that $G:=Aut(C)$ has an extra involution and $\u=(u_1, u_2)$    its dihedral
invariants.  Then,

a) $G\iso \bZ_3 \sem D_8$ if and only if  $(u_1, u_2 )=(0,0)$ or $(u_1, u_2 )=(6750, 450)$.

b) $G\iso GL_2(3) $ if and only if $(u_1, u_2) = ( -250, 50)$.

c) $G\iso D_{12}$ if and only if $\, \,  u_2^2 - 220 u_2 -16 u_1 +4500=0$, for $u_2 \neq 18, 140 + 60\sqrt{5},
50$.

d) $G\iso D_8$ if and only if         $\, \, 2 u_1^2-u_2^3=0$, for $u_2 \neq 2, 18, 0, 50, 450$. Cases $u_2 = 0,
450$ and $u=50 $ are reduced to cases a) and b) respectively.
\end{lemma}
  The mapping $\Phi: (u_1, u_2 ) \to (i_1, i_2, i_3)$,
gives a  birational parameterization of $\L_2$. The fibers of $\Phi$ of cardinality $>1$ correspond  to those
curves $C$ with $|Aut(C)| > 4$. Dihedral invariants  $u_1, u_2$ are given explicitly as rational functions of
$i_1, i_2, i_3$. The curve $Y^2=X^6-X$ is the only genus 2 curve (up to isomorphism) which has extra automorphisms
and is not in $\L_2$. The automorphism group in this case is $\Z_{10}$, see \cite{SV}. Thus, if $C\in \L_2$ we
determine $Aut(C)$ via Lemma 3.1., otherwise $C$ is isomorphic to $Y^2=X^6-X$ or $Aut(C)\iso \Z_2$.


The case $g=3$ is given as an application in \cite{GS2}. Let $C\in \L_3$ with equation as in Eq.~\ref{eq}.
Dihedral invariants are
$u_1= a_1^4+a_3^4,  u_2=(a_1^2+a_3^2)a_2,  u_3=2 a_1 a_3.$ The analogue of Lemma 3.1 is proved in \cite{GS2} for
$g=3$.

This technique can be used successfully for all $g$. We have implemented programs that determine $Aut(C)$ for
$C\in \L_g$ and for $g=2,3,4,5,6$. In order to compute the automorphism group of a curve $C \in \L_g$ we transform
this curve to its normal form (i.e., Eq.~\ref{eq}) and then compute its dihedral invariants. If these invariants
satisfy any locus $\L_G$ then the automorphism group is $G$, otherwise the automorphism group is $V_4$. The
following lemma determines a nice condition for $\bG$ to have at least two involutions.

\begin{lemma} For a curve $C\in \L_g$
the reduced automorphism group has at least two involutions if and only if
\begin{equation}\label{eqq}
(2^{g-1} u_1^2 + u_g^{g+1})\, (2^{g-1} u_1^2 - u_g^{g+1}) =0
\end{equation}
\end{lemma}
\begin{proof}
Let $C\in \L_g$. Then, there is an involution
 $z_1 \in \bar G$ which fixes no Weierstrass points
of $C$, see the proof of lemma 1 in \cite{GS2}. Thus, $z_1 (X)=-X$. Let $z_2\neq z_1$ be another involution in
$\bar G$. Since, $z_2 \neq z_1$ then $z_2(X)= \frac m X$, where $m^2=1$. Then, $V_4=\<z_1, z_2\> \emb \bar G$ and
$z_2$ or $z_1\, z_2$ is the transformation $X \to \frac 1 X$, say $z_2(X)=\frac 1 X$. If $g$ is odd we have $\P=\{
\pm \a_1, \pm \frac 1 {\a_1}, \dots , \pm \a_n, \pm \frac 1 {\a_n} \}$, where $n=[\frac {g+1} 2]$, otherwise $\P$
contains also two  points $\pm P$. Thus, $\pm P$ can be either fixed  or permuted by  $z_2(X)=\frac 1 X$. Hence,
 they are  $\pm 1$ or $\pm I$, where $I^2=1$.
The equation of $C$ is given by
\begin{eqnarray*}
& Y^2 =  &\prod_{i=1}^{n} (X^4-\lambda_i X^2 +1), \quad \textit{if g is odd}\\
& Y^2 =  & (X^2 \pm 1) \, \prod_{i=1}^{n} (X^4-\lambda_i X^2 +1), \quad \textit{if g is even}.
\end{eqnarray*}
Let $s:=\l_1 + \dots + \l_n$. If $g$ is odd then  $a_1=a_g=-s$.   Then, $u_1=2s^{g+1}$ and $u_g=2 s^2$ and they
satisfy Eq.~\ref{eqq}. If $z_2(X)=\frac 1 X$ fixes two points of $\P$ then one of the factors of the equation is
$X^4-1$. Then, $a_1= (-1)^{\frac 1 {g+1}}\, s$ and $a_g=(-1)^{\frac 1 {g+1}} \, s$. Hence, $a_g^{g+1}=
a_g^{g+1}=-s^{g+1} $ and $u_1=-2 s^{g+1}$, $u_g= -2 s^2$. Then, $2^{g-1} u_1^2 + u_g^{g+1}=0$.

If $g$ is even and $\{\pm 1\} \subset \P$  then
 $a_1= a_g=s+1$.  If $\{\pm I\} \subset \P$ then $a_1=a_g=1-s$.
In both cases  $2^{g-1} u_1^2 - u_g^{g+1}=0$. The converse goes similarly.

\end{proof}
\begin{remark}
If $2^{g-1} u_1^2 + u_g^{g+1}=0$, then one of the involutions $z_2$, $z_1 z_2$ of $\bar G$ lifts to an element  of
order 4 in $G$. If $ 2^{g-1} u_1^2 - u_g^{g+1} =0$ both of them lift to involutions in $G$.
\end{remark}
%
For $C \not \in \L_g$ we  check if $C$ has automorphisms of order $3 \leq N \leq 2(2g+1)$, see  Wiman \cite{Wi}.
The following lemma is a  consequence of \cite{Bu} and gives possible values for $N$.  We only sketch the proof.
\begin{lemma}\label{lemma3.3}
Let $C$ be a genus $g$ hyperelliptic curve with an automorphism of order $N>2$. Then  either $N=3, 4$ or one of
the following holds;

i)  $N|(2g+1)$ or $N|2g$ and $N<g$ (then $Aut(C)\iso \Z_{2N}$)

ii)   $N|2g$ and $N$ is an even number such that $ 6 \leq N \leq 2g-2$.

iii)  $N=4N'$ such that   $N'|g$ and $N'< g$.
\end{lemma}
%
\begin{proof}
 Let $C$ be a genus $g$ hyperelliptic curve with extra automorphisms such that
$C \not \in \L_g$. Then,  the automorphism group of $C$ is isomorphic to  one  of  the   following:
 $SL_2(3)$,  $SL_2(5)$,  $W_3$, $H_{N/2}$, $U_{N/2}$, $G_{N/2}$,
$\Z_{2N}$ where $N\, |\, 2g+1$ or $N\, | \,2g $ and $N < g$; see \cite{Bu} for definitions of these groups. All
other groups  listed in  Table 2  in \cite{Bu}  contain at least  two involutions,  hence  they correspond  to
curves in  $\L_g$. The only groups in the above list that might not contain an element of order 2, 3, or 4 are
$U_{N/2}$, $G_{N/2}$. The group $G_{N/2}$  (resp., $U_{N/2}$) has an element of order $N$ where $N$ is as above.

\end{proof}
To have a complete algorithm that works for any $g\geq 2$,  one needs to classify (up to isomorphism) curves of
genus $g$ which are not in the locus $\L_g$.
 In order to do this, we need invariants which
classify isomorphism classes of curves with an automorphism of order $N>2$.
 However, for small genus ad-hoc methods can be used to identify
such groups.

\section{Field of moduli}

In this section we introduce a method to compute the field of moduli of hyperelliptic curves with extra
automorphisms. Until recently this was an open problem even for $g=2$. Further, we state some open questions for
higher genus and prove Conjecture 1 for $\p \in \H_g$ such that the reduced automorphism group of $\p$ has at
least two involutions.

 Let $C$ be a genus $g$ hyperelliptic curve defined over $k$.
We can write the equation of $C$ as follows
$$Y^2=X(X-1)(X^{2g-1} + c_{2g-2} X^{2g-2} + \dots + c_1 X + c_0)$$
where the discriminant $\Delta$ of the right side is nonzero. Then, there is a map
\begin{eqnarray*}
& \Phi_1: \quad k^{2g-1}\setminus\{\Delta \neq 0\} \to \H_g\\
& \qquad (c_0, \dots , c_{2g-2}) \to \p=[C]
\end{eqnarray*}
of degree $d=4g(g+1)(2g+1)$. We denote by $J_{\Phi}$ the Jacobian matrix of a map $\Phi$. Then Conjecture 1 can be
stated as follows:

\medskip

\noindent {\bf Conjecture 2:} {\it For each $\p$  in the locus $ \det (J_{\Phi_1})=0$ such that $\p \in \H_g(L)$
there  exists a  representative $C$ of the isomorphism  class $\p$ which is defined over $L$. }

\medskip

For $g=2$ this conjecture is a theorem  as shown in \cite{CQ}. The main result in \cite{CQ} is to prove the case
when automorphism group is $V_4$. A method of Mestre is generalized which uses covariants of order 2 of binary
sextics and a result of Clebsch. Such a method probably could be generalized to higher genus as claimed by Mestre
\cite{Me} and Weber \cite{We}.
\begin{remark}
There is a mistake in the proof of Theorem 2 in \cite{CQ}. In other words, the proof is incorrect  when the
Clebsch invariant $C_{10}=0$. However, it can easily be fixed.  A correct version of the algorithm has been
implemented in Magma by P. van Wamelen.
\end{remark}
For $g=3$ the conjecture is proven by Gutierrez and this author for all points $\p$ with $|Aut(\p)|>4$, see
\cite{GS2}. The proof uses dihedral invariants of hyperelliptic curves. A generalization of the method used in
\cite{Me}, \cite{We} for $\p\in \H_3$ such that $Aut(\p)\iso V_4$ would complete the case $g=3$.

Next we focus on the locus $\L_g$.   Let   $C \in \L_g$. Then,
 $C$ can be written in the normal form as in
equation \ref{eq}. The map
\begin{eqnarray*}
& \Phi: \quad k^{g} \setminus \{\Delta \neq 0\} \to \L_g\\
& \qquad \qquad (a_1, \dots , a_g) \to (u_1, \dots , u_g)
\end{eqnarray*}
has  degree $d=2g+2$. We  ask a similar question as in Conjecture 2. Let $\p$ be in the locus $\,  \det
(J_{\Phi_1})=0\, $ such that $\p \in \H_g(L)$. Is  there   a  representative $C$ of the isomorphism  class $\p$
which is defined over $L$?

The determinant of the Jacobian matrix is
$$\det  (J_{\Phi} )=(2^{g-1} u_1^2 + u_g^{g+1})\, (2^{g-1} u_1^2 - u_g^{g+1}).$$
The locus $\det (J_{\Phi} )= 0 \,$   corresponds exactly to the hyperelliptic curves with  $V_4 \emb \bar G$ as
shown by Lemma 3.2.

\begin{theorem}
For each $\p$  in the locus $ \det (J_{\Phi})=0$ such that $\p \in \H_g(L)$ there  exists a  representative $C$ of
the isomorphism  class $\p$ which is defined over $L$.
 Moreover, the equation of $C$ over $L$
is given by
\begin{equation}
C: \quad Y^2=u_1 X^{2g+2} + u_1 X^{2g} + u_2 X^{2g-2} + \dots \pm  u_g X^2 +2,
\end{equation}
where the coefficient of $X^2$ is $u_g$ (resp., $-u_g$) when $2^{g-1} u_1^2-u_g^{g+1}=0$ (resp., $2^{g-1}
u_1^2+u_g^{g+1}=0$).
\end{theorem}

\begin{proof} Let $\p=(u_1, \dots , u_g)\in \L_g(L)$ such that $2^{g-1} u_1^2-u_g^{g+1}=0$.
All we need to show is that the dihedral invariants of $C$ satisfy the locus $\det (J_\Phi)=0$. By the appropriate
transformation $C$ can be written as
$$Y^2= X^{2g+2} + ( \frac {u_1} 2 )^{\frac 1 {g+1}} \cdot X^{2g} +
\sum_{i=1}^{g-1} \frac {u_{g+1-i}} {u_1} \cdot (\frac {u_1} 2 )^{ \frac {g+1-i} {g+1}} \cdot X^{2i}+1.$$
Then, its dihedral invariants are
$$u_1(C)= \frac {u_1} 2 + (\frac {u_g} {u_1} )^{g+1} \cdot ( \frac {u_1} 2)^g=
\frac { 2^{g-1} u_1^2 + u_g^{g+1} } {2^g u_1}, \quad u_g(C)=u_g. $$
Substituting $u_g^{g+1}=2^{g-1} u_1^2$ we get  $u_1(C)=u_1$. Thus, $C$ is in the isomorphism class determined by
$\p$ and defined over $L$.

Let $\p=(u_1, \dots , u_g)\in \L_g(L)$ such that $2^{g-1} u_1^2+u_g^{g+1}=0$. This case occurs only when $g$ is
odd, see the proof of Lemma 3.2. We transform $C$ as above and have $u_1(C)= u_1$ and $u_g(C)=-u_g$. They are the
other tuple $(u_1, ... , -u_g)$ which correspond to $\p$. This completes the proof.

\end{proof}
\noindent The following is a consequence of Lemma 3.2. and Theorem 4.1.
\begin{corollary}
Conjecture 1   holds for all $p\in \L_g$ such that the reduced automorphism group of $\p$ has at least two
involutions.
\end{corollary}

\section{Closing remarks}
Conjecture 1 was stated for the first time during a talk of the author in {\sc ANTS V}, see \cite{Sh1}. It can be
generalized to $\M_g$ instead of $\H_g$. However, little is known about the loci $\M_G$ (i.e., locus of curves in
$\M_G$ with full automorphism group $G$). In \cite{MS} we introduce an algorithm that would classify such groups
$G$ for all $g$ and give a complete  list of ``large'' groups for $g\leq 10$. However, finding invariants that
classify curves with automorphism group $G$ is not an easy task, since the equations describing non-hyperelliptic
curves are more complicated then the hyperelliptic case. A more theoretical approach on singular points of $\M_g$
probably would produce better results on Conjecture 1. At this time we are not aware of any such results.

Our approach would work (with necessary adjustments) even in positive characteristic. However, the goal of this
note was to introduce such method rather than
 explore it to the full extent.

Computationally, dihedral  invariants  give an efficient way  of determining a point of the moduli space $\L_g$.
Using such invariants in positive characteristic
 could have applications in the arithmetic of hyperelliptic curves,
 including  cryptography.

\section*{Acknowledgments}
This paper was written during a  visit  at the University of Florida. I want  to thank  the Department  of
Mathematics at  the University  of Florida for their hospitality.
%


\end{document}